\documentclass[titlepage,12pt]{article}

\usepackage{graphicx}
\usepackage{amssymb}

\newcommand{\R}{I\!\!R}

\begin{document}

\baselineskip=18pt

\vspace*{1cm}

\begin{center}
{\Large{\bf The Aharonov-Bohm effect:\\ Mathematical Aspects of
the Quantum Flow}}
\end{center}

\vspace{0.3cm}

\begin{center}
{\large Luis Fernando Mello}
\end{center}
\begin{center}
{\em Instituto de Ci\^encias Exatas, Universidade Federal de
Itajub\'a \\CEP 37.500-903, Itajub\'a, MG, Brazil
\\}e--mail:lfmelo@unifei.edu.br
\end{center}

\vspace{0.2cm}

\begin{center}
{\large Yuri C\^andido Ribeiro}
\end{center}
\begin{center}
{\em Instituto de Sistemas El\'etricos e Energia, Universidade
Federal de Itajub\'a
\\CEP 37.500-903, Itajub\'a, MG, Brazil
\\}e--mail:yuricandido@yahoo.com.br
\end{center}

\vspace{0.05cm}

\begin{center}
{\bf Abstract}
\end{center}

This paper addresses the scattering of a beam of charged particles
by an infinitely long magnetic string in the context of the
hydrodynamical approach to quantum mechanics. The scattering is
qualitatively analyzed by two approaches. In the first approach,
the quantum flow is studied via a one-parameter family of complex
potentials. In the second approach, the qualitative theory of
planar differential equations is used to obtain a one-parameter
family of Hamiltonian functions which determine the phase
portraits of the systems.

\vspace{0.3cm}



\newpage
\baselineskip=20pt

\section{\bf Introduction}\label{intro}

Consider a long solenoid of small transverse cross section endowed
with a magnetic flux. The limit configuration when the cross
section becomes vanishingly small while the magnetic flux enclosed
is kept constant is called a {\it magnetic string}. The magnetic
field vanishes everywhere except inside the magnetic string and,
as there is no Lorentz force, charged particles around the string
are not affected by it, according to the classical mechanics point
of view. Nevertheless, there is a nontrivial quantum scattering
due to the Aharonov-Bohm effect. See Ref. \cite{m} and \cite{op}.

The nontrivial scattering mentioned above can be studied from the
following family of vector fields (the probability current
density) (see \cite{op})
\begin{equation}
\mathcal J \left( x,y,\delta \right)= \frac{\hbar k}{M} \left( -
1+ \displaystyle \frac{\delta}{k} \:
\frac{y}{x^2+y^2},-\displaystyle \frac{\delta}{k} \:
\frac{x}{x^2+y^2} \right), \label{current}
\end{equation}
or equivalently by the following family of planar ordinary
differential equations
\begin{equation}
\left\{ \begin{array}{l}
x'= \displaystyle \frac{\hbar k}{M} \left(-1 + \displaystyle \frac{\delta}{k} \frac{y}{x^2+y^2}\right) \\
\\
y'= - \displaystyle \frac{\hbar k}{M} \left( \displaystyle
\frac{\delta}{k} \frac{x}{x^2+y^2} \right)
\end{array} \right. \label{equat1}
\end{equation}
where $M$ is the mass of the particle, $0 \leq k < \infty $ is
associated to the energy $ \hbar^2 k^2 / 2 M$ for a stationary
state and $0 \leq \delta \leq \frac{1}{2}$ is the flux parameter
(see section \ref{tube}).

In \cite{m} the numerical solution of the system (\ref{equat1}) is
plotted, showing the main features of the quantum flow near the
magnetic string. When $\delta = 0$ the streamlines are parallel to
the $x$-axis which is identical to the classical flow (see Fig.
\ref{paralelo}). However, when $\delta > 0$ the topology of the
streamlines changes near the origin, where open lines become
cycles. As Fig. \ref{estag1} shows, the quantum and the classical
flow differ radically from each other.

This paper provides a new analysis of the scattering problem near
the magnetic string, using complex potentials and the qualitative
theory of planar differential equations. In Section \ref{tube} we
introduce some notation and definitions. The family of velocity
fields (\ref{current}) is analyzed via complex potentials in
Section \ref{comppot}, showing that the quantum flow of
(\ref{current}) is determined by the streamlines of a family of
complex potentials. In Section \ref{hamil}, the system
(\ref{equat1}) is studied using the qualitative theory of planar
differential equations. We show that the family of differential
equations (\ref{equat1}) is a family of planar Hamiltonian systems
whose phase portraits are obtained by the level curves of the
Hamiltonian functions. Conclusions are presented in Section
\ref{conclusion}, and a review on complex potentials is given in
the Appendix.

\section{The magnetic string}\label{tube}

Let us assume an infinity magnetic string carrying magnetic flux
$\Phi$ coinciding with the $z$ axis. The vector potential is then
given by
\begin{equation}
{\bf A}= \frac{\Phi}{2 \pi r} \: {\it e}_\theta, \label{magnetic}
\end{equation}
where ${\it e}_\theta = (-\sin \theta, \cos \theta)$, $r$ and
$\theta$ are the polar coordinates in the plane $xy$. The flux
parameter $\delta$ is defined by
\begin{equation}
\delta = \frac{e \Phi}{2 \pi \hbar c}. \label{delta}
\end{equation}

The evolution of the state $\Psi$ of a particle of charge $e$ and
mass $M$ interacting with an electromagnetic field of scalar
potential $\varphi$ and vector potential ${\bf A}$ is given by the
Schr\"{o}dinger equation
\begin{equation}
i \: \hbar \: \frac{\partial \Psi}{\partial t}= \hat{H} \: \Psi,
\label{schrodinger}
\end{equation}
where the Hamiltonian operator $\hat{H}$ is
\begin{equation}
\hat{H}= \frac{1}{2M} \left[ -i \: \hbar \: \nabla - \frac{e}{c}
\: {\bf A} \right]^2 + e \: \varphi. \label{hamiltonoperator}
\end{equation}

The hydrodynamical approach (see \cite{op}) is obtained by writing
the wave function in the form
\[
\Psi = \rho \: e^{i\chi},
\]
and separating the real and imaginary parts of the Schr\"{o}dinger
equation. The real part leads to
\begin{equation}
M \left(\frac{\partial {\bf v}}{\partial t}+ ({\bf v} \cdot
\nabla){\bf v} \right) = e \: {\bf E}+\frac{e}{c} \: {\bf v}
\times {\bf B} - \nabla V, \label{real}
\end{equation}
where
\begin{equation}
{\bf v}= \frac{1}{M} \left(\hbar \nabla \chi - \frac{e}{c} \: {\bf
A}\right), \label{veloc}
\end{equation}
and
\begin{equation}
V= -\frac{\hbar}{2M}\frac{\Delta \rho}{\rho}, \label{quantumpot}
\end{equation}
is the {\it quantum potential}. The imaginary part yields
\begin{equation}
\frac{\partial \rho^2}{\partial t} + \nabla \cdot (\rho^2 {\bf
v})=0, \label{imaginaria}
\end{equation}
where $\rho^2$ is the {\it probability density} and $\mathcal J=
\rho^2 \: {\bf v}$ is the {\it probability current}.

In the scattering problem near the magnetic string, Eq.
(\ref{real}) and (\ref{imaginaria}) have the forms
\begin{equation}
M \left(\frac{\partial {\bf v}}{\partial t}+ ({\bf v} \cdot
\nabla){\bf v} \right) = - \nabla V, \label{real1}
\end{equation}
\begin{equation}
\nabla \cdot \mathcal J =0, \label{imaginaria1}
\end{equation}
since ${\bf E}=0$ and ${\bf B}=0$, whereas the probability current
associated with the Aharonov-Bohm wave function is given by Eq.
(\ref{current}) (see \cite{op}).

\section{Complex potentials}\label{comppot}

\newtheorem{teo}{Theorem}[section]
\newtheorem{lema}[teo]{Lemma}
\newtheorem{prop}[teo]{Proposition}
\newtheorem{cor}[teo]{Corollary}

We begin with the family of complex analytic functions
\begin{equation}
F \left( z,\delta \right)= \frac{\hbar k}{M} \left( -z +
\displaystyle i \: \frac{\delta}{k} \: \log z \right),
\label{potent}
\end{equation}
where $z = x + i y$. It follows that
\begin{equation}
F^{\prime} \left( z,\delta \right)= \frac{\hbar k}{M} \left( -1 +
i \: \frac{\delta}{k} \: \frac{1}{z} \right), \label{F'}
\end{equation}
which can be written as
\begin{equation}
F^{\prime} \left( z,\delta \right)= \frac{\hbar k}{M} \left[
\left(-1 + \frac{\delta}{k} \: \frac{y}{x^2+y^2} \right) - i \:
\left(- \frac{\delta}{k} \: \frac{x}{x^2+y^2} \right) \right].
\label{F'1}
\end{equation}
From Eqs. (\ref{def F'}) and (\ref{F'1}) it follows that
\begin{equation}
u \left( x,y,\delta \right)= \frac{\hbar k}{M} \left( -1 +
\frac{\delta}{k} \:\frac{y}{x^2+y^2} \right)
\end{equation}
and
\begin{equation}
v \left( x,y,\delta \right)= - \frac{\hbar k}{M} \left(
\frac{\delta}{k} \:\frac{x}{x^2+y^2} \right),
\end{equation}
are the components of probability current (\ref{current}).
Therefore we have proved the following lemma.

\begin{lema}
The family of complex analytic functions (\ref{potent}) is the
family of complex potentials of the probability current
(\ref{current}).

\label{lema1}
\end{lema}

The complex potential $F$ in Eq. (\ref{potent}) is the sum of the
{\it ``classical" complex potential} $F_1$ and the {\it ``quantum"
complex potential} $F_2$,
\begin{equation}
F_1(z) = - \frac{\hbar k}{M} \: z
\end{equation}
and
\begin{equation}
F_2 \left( z,\delta \right) = i \: \frac{\hbar \delta}{M} \: \log
z.
\end{equation}

\begin{figure}[!h]
\centerline{
\includegraphics[width=8cm]{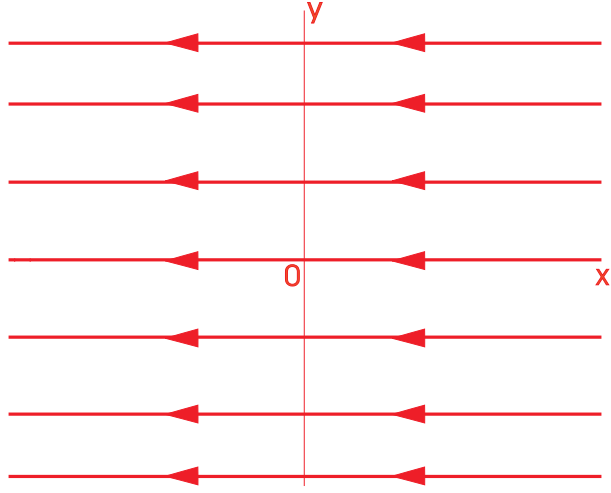}}
\caption{{\small Flow of (\ref{current}) when $\delta = 0$}.}

\label{paralelo}
\end{figure}

The stream function of the complex potential $F_1(z)= - (\hbar k /
M) \: z = (\hbar k / M) \: (-x - i y)$ is $\psi_1(x,y)= - (\hbar k
/ M) \: y$. Thus the streamlines are the horizontal lines $y = \:
c_0$ and the velocity at any point is $F'_1(z)= - \hbar k / M$.
The flow is uniform to the left and is called {\it parallel}. When
$\delta = 0$ we have $F = F_1$, and the flow is illustrated in
Fig. \ref{paralelo}.

\begin{figure}[!h]
\centerline{
\includegraphics[width=8cm]{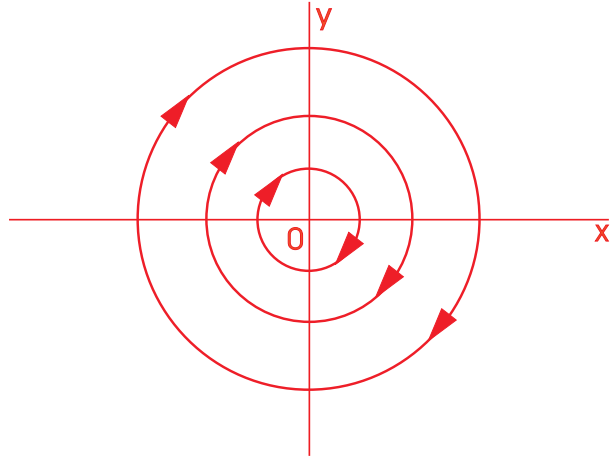}}
\caption{{\small Streamlines of quantum complex potential $F_2
$}.}

\label{rota}
\end{figure}

In polar coordinates the complex potential $F_2$ has the form
\begin{equation}
F_2 \left( r,\theta,\delta \right) = i \: \frac{\hbar \delta}{M}
\: \log (r {\rm e}^{\, i \, \theta}) = -\frac{\hbar \delta}{M} \:
\theta + i \: \frac{\hbar \delta}{M} \: \log r. \label{def F21}
\end{equation}
The stream function is
\begin{equation}
\psi_2 \left( r,\theta,\delta \right) = \frac{\hbar \delta}{M} \:
\log r
\end{equation}
and the streamlines are given by $\psi_2 \left( r,\theta,\delta
\right) = (\hbar \delta / M) \: \log r = \: c_0$, which are
circles centered at the origin as represented in Fig. \ref{rota}.
This flow is a {\it pure rotation} with a {\it vortex} at the
origin.

This analysis shows that the flow of the complex potential $F$ in
Eq. (\ref{potent}) is the superposition of the parallel and the
pure rotation flows. A new stagnation point appears to solve the
compatibility problem of the two flows. It follows from Eq.
(\ref{F'})
\[
F^{\prime}\left( z,\delta \right) =  -\frac{\hbar k}{M}\:\left( 1
- \frac{z_0}{z} \right),
\]
\noindent where $z_0 = i \: (\delta / k)$. Thus $F^{\prime}\left(
z,\delta \right) = 0$ if and only if $z = z_0$.

\begin{figure}[!h]
\centerline{
\includegraphics[width=8cm]{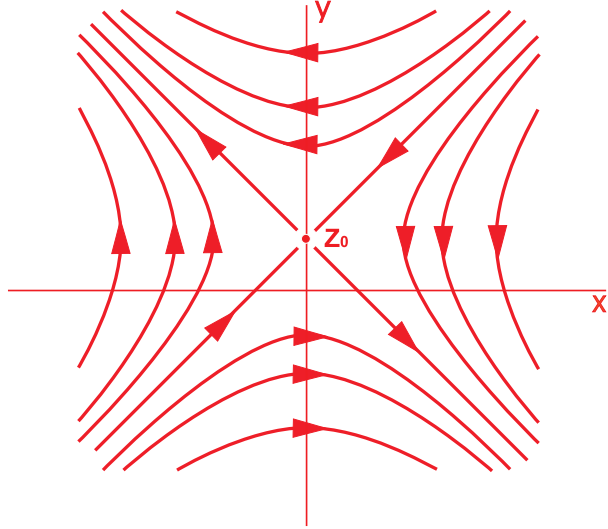}}
\caption{{\small Streamlines of the complex potential $F_3$}.}

\label{hipe}
\end{figure}

The Taylor expansion of $F'\left( z,\delta \right)$ at $z_0$ is
given by
\begin{equation}
F^{\prime} \left( z,\delta \right) = - \frac{\hbar k}{M} \:
\frac{z-z_0}{z_0} + O(|z|^2) = i \: \frac{\hbar k^2}{\delta M} \:
(z-z_0) + O(|z|^2).
\end{equation}
Therefore near of the stagnation point $z_0$ the complex potential
can be written as
\begin{equation}
F_3 \left( z,\delta \right) = i\: \frac{\hbar k^2}{2 \delta M}
 \: (z-z_0)^2.
\end{equation}
Take $z= z - z_0$ for the study of the complex potential $F_3$.
Thus
\begin{equation}
F_3 \left( z,\delta \right) = i \: \frac{\hbar k^2}{2 \delta M} \:
(x+i\:y)^2, \label{def F31}
\end{equation}
and the stream function is given by
\begin{equation}
\psi_3 \left( x,y,\delta \right) = \frac{\hbar k^2}{2 M \delta}
 \: (x^2-y^2). \label{parteimag F3}
\end{equation}
The streamlines are the level curves
\[
\frac{\hbar k^2}{2 M \delta}  \: (x^2-y^2) = \: c_0.
\]
By a translation, the flow near $z_0 = i \: (\delta /k)$ is
determined, as illustrated in Fig. \ref{hipe}. We have proved the
following lemma.

\begin{lema}
The flow of the complex potential $F$ in Eq. (\ref{potent}) is the
superposition of the parallel and pure rotation flows, as
represented in Fig. \ref{paralelo} and \ref{rota} respectively.
The complex potential $F$ has an stagnation point at $z_0 = i \:
(\delta / k)$ for all $\delta >0$. The flow of $F$ near of $z_0$
is represented in Fig. \ref{hipe}.

\label{lema2}
\end{lema}

One sees from Eq. (\ref{F'}) that
\[
\lim_{|z| \rightarrow  \infty} F'\left( z,\delta \right) = -
\frac{\hbar k}{M} = F'_1 (z).
\]
Thus as the distance from the origin increases the quantum flow
approaches that of a plane wave, as expected.

It should be noted that the stream function is given by
\begin{equation}
\psi \left( x,y,\delta \right) = \psi_1 (x,y) + \psi_2 \left(
x,y,\delta \right) = \frac{\hbar k}{M} \left( - y +
\frac{\delta}{k} \: \log \left(\sqrt {x^2+y^2} \right) \right).
\label{stream}
\end{equation}
Thus $\psi (-x,y,\delta) = \psi (x,y,\delta)$ and therefore the
quantum flow is symmetric with respect to $y$-axis.

\begin{figure}[!h]
\centerline{
\includegraphics[width=8cm]{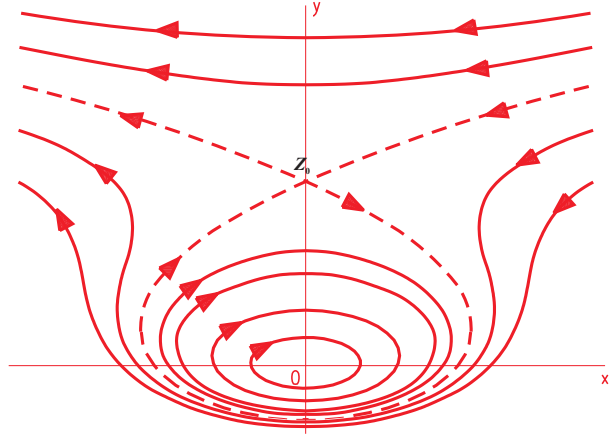}}
\caption{{\small  Quantum flow of (\ref{current}) near the
magnetic string ($\delta
> 0 $)}.}

\label{estag1}
\end{figure}

The separatrices of the stagnation point $z_0= i\: (\delta / k) $,
for $\delta > 0 $, represented in Fig. \ref{estag1} by dashed
lines, are subsets of the level curve
\[
L (\delta)= \left \{ \,(x,y) \in \R^2 \bigg | \: \psi ( x,y,
\delta) = \displaystyle \frac{\hbar \delta}{M} \: \left(\log
\displaystyle \frac{\delta}{k} \: -1 \right),\delta > 0 \right\}.
\]
There is a homoclinic streamline of the stagnation point $z_0$
bounding a neighborhood filled with cycles surrounding the vortex
$(0,0)$. This neighborhood shrinks to $(0,0)$ as $\delta \to 0^+
$.

The following theorem is an immediate consequence of our analysis.
\begin{teo}
The family of complex analytic functions $F \left( z,\delta
\right)$ in (\ref{potent}) is the family of quantum complex
potentials of the quantum flow of (\ref{current}). The quantum
flow of (\ref{current}) are represented in Fig. \ref{paralelo} and
Fig. \ref{estag1} for $\delta = 0$ and $\delta
> 0$, respectively.

\label{teopot}

\end{teo}
The circulation of $\mathcal J$ is given by
\begin{equation}
\Gamma_{\mathcal J} = \oint_C \mathcal J \cdot t \: ds = \oint_C
F'(z) \: dz = \oint_C \left[ \frac{\hbar k}{M} \left( -1 + i \:
\frac{\delta}{k} \frac{1}{z} \right) \right] \: dz = -2 \pi
\frac{\hbar \delta}{M}, \label{circul}
\end{equation}
where $C$ is a closed curve encircling the origin. From Eqs.
(\ref{delta}) and (\ref{circul}), this circulation can be
rewritten as
\begin{equation}
\Gamma_{\mathcal J} = - \frac{e}{cM} \Phi. \label{circul1}
\end{equation}

Therefore from the mathematical point of view the quantum flow
associated to the Aharonov-Bohm effect is a parallel flow at
infinity with a vortex at the origin, where the circulation around
the vortex is due to nonzero magnetic flux $\Phi$.

\section{Planar Hamiltonian differential equations}\label{hamil}

The planar differential equations of the form
\begin{equation}
\left\{ \begin{array}{l}
x'= \displaystyle \frac{\partial H}{\partial y}(x,y)\\
\\
y'=- \displaystyle \frac{\partial H}{\partial x}(x,y)
\end{array} \right.
\label{equat3}
\end{equation}
are called {\it Hamiltonian differential equations}, where $H$ is
a $C^2$ real function called a {\it Hamiltonian function}. The
Hamiltonian differential equations generally arise in mechanical
systems without friction, and the total energy of the system can
often be taken as the Hamiltonian function. See \cite{jl} and
references therein for more details about mathematical aspects of
Hamiltonian differential equations.

It is clear from the special form of Eq. (\ref{equat3}) that the
Hamiltonian $H$ is a first integral (conservation of energy). Thus
the phase portrait of Eq. (\ref{equat3}) is obtained by the level
curves of $H$.

The differential equations (\ref{equat1}) are Hamiltonian
differential equations with Hamiltonian
\begin{equation}
H ( x,y,\delta)= \displaystyle \frac{\hbar k}{M} \left(
\frac{\delta}{k} \: \log \left(\sqrt {x^2+y^2} \right)-y \right).
\label{equat4}
\end{equation}
Observing Eqs. (\ref{stream}) and (\ref{equat4}), one sees that
$H(x,y,\delta)= \psi (x,y,\delta)$. We have the following analog
of the Theorem \ref{teopot}.

\begin{teo}
The planar differential equations (\ref{equat1}) are Hamiltonian
differential equations with the Hamiltonian $H$ given by
(\ref{equat4}). The phase portraits of the differential equations
(\ref{equat1}) are illustrated in Fig. \ref{paralelo} and Fig.
\ref{estag1} for $\delta = 0$ and $\delta > 0$, respectively.

\label{teohamilt}

\end{teo}

Therefore, in the context of qualitative theory of differential
equations, (\ref{current}) is a Hamiltonian vector field where the
stream function plays the role of the Hamiltonian.

\section{Concluding remarks}\label{conclusion}

In sections \ref{comppot} and \ref{hamil} we have analyzed, by two
different points of view, the mathematical aspects of the quantum
flow associated with the Aharonov-Bohm effect. In section
\ref{comppot} we use the complex potentials theory to obtain the
features of the velocity fields (\ref{current}). The main goal was
to obtain the family of complex potentials (\ref{potent}).

In section \ref{hamil}, we used the qualitative theory of planar
differential equations to obtain the phase portraits of the planar
differential equations (\ref{equat1}), showing that the
differential equations (\ref{equat1}) are planar Hamiltonian
differential equations.

\section*{Appendix: A review on complex potentials}

\renewcommand{\theequation}{A.\arabic{equation}}

Complex potentials play an important role in hydrodynamics
\cite{mh}. Here the motion of the fluid is assumed to be
two-dimensional, i.e., the same in all planes parallel to the
$xy$-plane, and is independent of time. Thus, it is sufficient to
consider the motion in the $xy$-plane. We consider only
irrotational flows, and we also assume that the fluid is
incompressible and free from viscosity.

Let $q(x,y)=(u(x,y),v(x,y))$ be the velocity of a particle of the
fluid at any point $(x,y)$. Here the functions $u(x,y)$ and
$v(x,y)$ are of class $C^1$. From the above hypotheses
\[
\oint_C q \cdot n \:\: ds = 0,
\]
where $C$ is a positively oriented simple closed contour lying in
a simply connected domain, $n$ is the normal vector to $C$ and
$ds$ is the arc length of $C$. From Gauss Theorem
\begin{equation}
{\nabla \cdot q(x,y)}=\frac{\partial u}{\partial x} (x,y) +
\frac{\partial v}{\partial y}(x,y)=0. \label{div nulo}
\end{equation}

As the flow is irrotational one finds
\begin{equation}
\oint_C q \cdot t \:\: ds = 0, \label{circulacao}
\end{equation}
where $t$ is the tangent vector to $C$. From Green Theorem
\begin{equation}
\nabla \times q(x,y) = \frac{\partial v}{\partial x} (x,y) -
\frac{\partial u}{\partial y}(x,y)=0. \label{circ nula}
\end{equation}

As Eqs. (\ref{div nulo}) and (\ref{circ nula}) are the
Cauchy-Riemann equations for the functions $u$ and $-v$,
\begin{equation}
u(x,y)-iv(x,y) \label{u+iv}
\end{equation}
is a complex analytic function.

It follows from Eq. (\ref{circulacao}) that $u \: dx + v \: dy$ is
an exact differential. Therefore there is a function $\phi (x,y)$,
the {\it velocity potential}, such that
\begin{equation}
q(x,y) = \nabla \phi (x,y). \label{def phi}
\end{equation}

From (\ref{def phi}) and (\ref{div nulo}) it follows that
\begin{equation}
\Delta \phi (x,y)=\frac{\partial^2 \phi}{\partial x^2} (x,y) +
\frac{\partial^2 \phi}{\partial y^2}(x,y)=0. \label{laplace}
\end{equation}
In other words, $\phi$ is a harmonic function. Let $\psi(x,y)$ be
a harmonic conjugate of $\phi(x,y)$. Thus
\begin{equation}
F = \phi  + i\psi \label{potdefinition}
\end{equation}
is a complex analytic function, called {\it complex potential} of
the flow.

As the level curves of $\psi$ and $\phi$ are orthogonal at points
where the velocity vector is not the zero vector, the velocity
vector is tangent to the level curves of $\psi$. Thus the function
$\psi$ characterizes the flow in a region. This function $\psi$ is
called {\it stream function} and their level curves are called the
{\it streamlines} of the flow. A point $z$ where $F^{\prime} (z) =
0$ is called an {\it stagnation point}.

From Eq. (\ref{potdefinition}) and (\ref{u+iv})
\begin{equation}
F^{\prime} = u - i v = \bar q. \label{def F'}
\end{equation}

\vspace{0.1cm} \noindent {\bf Acknowledgments.} The first author
developed this work under the project FAPEMIG CEX-1597/2005. The
authors are grateful to E. S. Moreira Jr. and N. Manzanares Filho
for reviewing the manuscript and for clarifying discussions.

\end{document}